\documentclass{amsart}
\usepackage{graphicx,amssymb,stmaryrd,amsfonts,epsfig,amsthm,a4,amsmath,mathrsfs,url}
\usepackage{vmargin,amscd}
\usepackage{eufrak}
\usepackage[latin1]{inputenc}
\psfigdriver{dvips}

\newtheorem{thm}{Theorem}[section]
\newtheorem{cor}[thm]{Corollary}
\newtheorem{clai}[thm]{Claim}
\newtheorem{lem}[thm]{Lemma}

\newtheorem{prop}[thm]{Proposition}
\theoremstyle{definition}
\newtheorem{defn}[thm]{Definition}
\theoremstyle{remark}
\newtheorem{rem}[thm]{Remark}
\newtheorem{exe}[thm]{Example}
\numberwithin{equation}{section}

\newcommand{\eps}{\varepsilon}

\newcommand{\A}{\mathcal{A}}

\newcommand{\Z}{\mathbf{Z}}

\newcommand{\N}{\mathbf{N}}
\newcommand{\R}{\mathbf{R}}

\newcommand{\tpr}{\begin{tiny}\noindent Proof:}

\newcommand{\bpr}{\noindent \textbf{Proof}: ~}
\newcommand{\epr}{~$\blacksquare$}
\begin{document}
\title{Asymptotic isoperimetry of balls in metric measure spaces}%
\author{Romain Tessera }
\date{\today}

\newpage

\begin{abstract}
In this paper, we study the asymptotic behavior of the volume of
spheres in metric measure spaces. We first introduce a general
setting adapted to the study of asymptotic isoperimetry in a
general class of metric measure spaces. Let $\A$ be a family of
subsets of a metric measure space $(X,d,\mu)$, with finite,
unbounded volume. For $t>0$, we define:
$$I^{\downarrow}_{\A}(t)=\inf_{A\in \A, \mu(A)\geq t} \mu(\partial
A).$$ We say that $\A$ is asymptotically isoperimetric if $\forall
t>0$:
$$I_{\A}^{\downarrow}(t)\leq CI(Ct),$$
where $I$ is the profile of $X$. We show that there exist graphs
with uniform polynomial growth whose balls are not asymptotically
isoperimetric and we discuss the stability of related properties
under quasi-isometries. Finally, we study the asymptotically
isoperimetric properties of connected subsets in a metric measure
space. In particular, we build graphs with uniform polynomial
growth whose connected subsets are not asymptotically
isoperimetric.
\end{abstract}
\maketitle

\section{Introduction}
\subsection{Asymptotic isoperimetry in metric measure spaces.}
\subsubsection{Boundary of a subset and isoperimetric
profile.\\}\label{bord}

Let $(X,d,\mu)$ be a metric measure space. Let us denote by
$B(x,r)$ the open ball of center $x$ and radius $r$. We suppose
that the measure $\mu$ is Borel and $\sigma$-finite. For any
measurable subset $A$ of $X$, any $h>0$, write:
$$A_h=\{x\in X, d(x,A)\leq h\},$$
and:
$$\partial_h A=A_{h}\cap (A^{c})_h.$$

Let us call $\partial_h A$ the $h$-boundary of $A$, and
$\partial_h B(x,r)$ the $h$-sphere of center $x$ and radius $r$.

\begin{defn}
Let us call the $h$-profile the nondecreasing function defined on
$\R_+$ by:
$$I_h(t)=\inf_{\mu(A)\geq t}\mu(\partial_h A),$$
where $A$ ranges over all $\mu$-measurable subsets of $X$ with
finite measure.
\end{defn}
This definition of large-scale boundary has the following
advantage: under some weak properties on the metric measure space
$X$, we will see in Section \ref{loc} that in some sense, the
boundary of a subset $A\subset X$ has a thickness ``uniformly
comparable to $h$". This will be play a crucial role in the proof
of the invariance of ``asymptotic isoperimetric properties" under
large-scale equivalence.

We could also define the boundary of a subset $A$ by $\partial_h
A=A_h\smallsetminus A$. But with this definition, the thickness of
the boundary may have uncontrollable ``fluctuations". Indeed,
consider a Riemannian manifold $X$ with boundary and take a subset
$A$ which is at distance $\eps>0$ from the boundary of $X$. Assume
that $\eps\ll h$. With the latter definition, the thickness of
$\partial_h A$ may vary between $h$ and $\eps$ and its volume may
strongly depend on $\eps$, even if $X$ has bounded geometry. Note
that this problem disappears with our definition since every point
of the boundary of $A$ is close to a ball of radius $h/2$ included
in $\partial_h A$ (see Section \ref{loc} for precise statements).

\subsubsection{Lower/upper profile restricted to a family of subsets.\\}

Let $(X,d,\mu)$ be a metric measure space. In order to study
isoperimetric properties of a family of subsets of $X$ with
finite, unbounded volume, it is useful to introduce the following
notions:
\begin{defn}
Let $\A$ be a family of subsets of $X$ with finite, unbounded
volume. We call lower (resp. upper) $h$-profile restricted to $\A$
the nondecreasing function $I^{\downarrow}_{h,\A}$ defined by:
$$I^{\downarrow}_{h,\A}(t)=\inf_{\mu(A)\geq t, A\in \A}\mu(\partial_h A)$$
(resp. $I^{\uparrow}_{h,\A}(t)=\sup_{\mu(A)\leq t, A\in
\A}\mu(\partial_h A)$).
\end{defn}

\begin{defn}
Consider two monotone functions $f$ and $g$: $\R_+\rightarrow
\R_+$. Say that $f\approx g$ if there exist some constants $C_{i}$
such that $C_{1}f(C_{2}t)\leq g(t)\leq C_{3}f(C_{4}t)$ for all
$t\in\R_+$.
\end{defn}
The asymptotic behavior of a monotone function
$\R_+\rightarrow\R_+$ may be defined as its equivalence class
modulo $\approx$.

We get a natural order relation on the set of equivalence classes
modulo $\approx$ of monotone functions defined on $\R_+$ by
setting:
$$(f\preceq g)\Leftrightarrow (\exists C_1,C_2>0,\forall t>0,\quad f(t)\leq C_1g(C_2t)).$$

\ We say that the family $\A$ is {\bf asymptotically
isoperimetric} (resp. {\bf strongly asymptotically isoperimetric})
if for all $A\in \A$
$$I^{\downarrow}_{h,\A}\preceq I_h$$
(resp. $I^{\uparrow}_{h,\A}\preceq I_h$).

\begin{rem}
Note that asymptotically isoperimetric means that for any $t$ we
can always choose an optimal set among those of $\A$ whose measure
is larger than $t$ whereas strongly asymptotically isoperimetric
means that every set of $\A$ is optimal (but the family
$(\mu(A))_{A\in \A}$ may be lacunar). In almost all cases we will
consider, the family $(\mu(A))_{A\in \A}$ will not be lacunar, and
strong asymptotic isoperimetry will imply asymptotic isoperimetry.
\end{rem}

\subsubsection{Large scale study.\\}

Let us recall the definition of a quasi-isometry (which is also
sometimes called rough isometry).

\begin{defn}
Let $(X,d)$ and $(X',d')$ two metric spaces. One says that $X$ and
$X'$ are quasi-isometric if there is a function $f$ from $X$ to
$X'$ with the following properties.
\begin{itemize}
\item[(a)]there exists $C_1>0$ such that
$\left[f(X)\right]_{C_1}=X'.$

\item[(b)] there exists $C_2\geq 1$ such that, for all $x,y\in X$,
$$C_2^{-1}d(x,y)-C_2\leq d'(f(x),f(y))\leq C_2d(x,y)+C_2.$$
\end{itemize}

\end{defn}

\begin{exe}
Let $G$ be a finitely generated group and let $S_1$ and $S_2$ two
finite symmetric generating sets of $G$. Then it is very simple to
see that the identity map: $G\rightarrow G$ induces a
quasi-isometry between the Cayley graphs $(G,S_1)$ and $(G,S_2)$.
At the beginning of the 80's, M. Gromov (see \cite{G}) initiated
the study of finitely generated groups up to quasi-isometry.
\end{exe}

\begin{exe}
The universal cover of a compact Riemannian manifold is
quasi-isometric to every Cayley graph of the covering group (see
\cite{G} and \cite{Sal}).
\end{exe}

Note that the notion of quasi-isometry is purely metric. So, when
we look for quasi-isometry invariant properties of a metric
measure space like, for instance, volume growth, we are led to
assume some uniformity properties on the volume of balls. This is
the reason why, for instance, this notion is well adapted to
geometric group theory. But since we want to deal with more
general spaces, we will define a more restrictive class of maps.
Those maps will be asked to preserve locally the volume of balls.
On the other hand, we want local properties to be stable under
bounded fluctuations of the metric. Precisely, let $(X,d,\mu)$ be
a metric measure space and let $d'$ be another metric on $X$ such
that $d/d'$ and $d'/d$ are bounded. The following definition (see
\cite{CS}) prevents wild changes of the volume of balls with
bounded radii under the identity map between $(X,d,\mu)$ and
$(X,d',\mu)$.

\begin{defn}
Let us say that $(X,d,\mu)$ is {\bf doubling at fixed radius}, or
has property $(DV)_{loc}$ if for all $r>0$, there exists $C_r>0$
such that, for all $x\in X$:
$$\mu(B(x,2r))\leq C_r\mu(B(x,r)).$$
\end{defn}

\begin{rem}
Note that Property $(DV)_{loc}$ is local in $r$ but uniform in
$x$.
\end{rem}

\begin{exe}
Bounded degree graphs or Riemanniann manifolds with Ricci
curvature bounded from below satisfy $(DV)_{loc}$.
\end{exe}

The following notion was introduced by Kanai \cite{Kan} (see also
\cite{CS}).

\begin{defn}
Let $(X,d,\mu)$ and $(X',d',\mu')$ two metric measure spaces with
property $(DV)_{loc}$. Let us say that $X$ and $X'$ are large
scale equivalent (we can easily check that it is an equivalence
relation) if there is a function $f$ from $X$ to $X'$ with the
following properties: there exist some constants $C_1>0$, $C_2\geq
1$, $C_3\geq 1$ such that:
\begin{itemize}

\item[(a)] $f$ is a quasi-isometry of constants $C_1$ and $C_2$;

\item[(b)]for all $x\in X$:
$$C_3^{-1}\mu(B(x,1))\leq \mu'(B(f(x),1))\leq C_3\mu(B(x,1)).$$

\end{itemize}

\end{defn}
Focusing our attention on balls of radius $1$ may not seem very
natural. Nevertheless, this is not a serious issue since property
$(DV)_{loc}$ allows to make no distinction between balls of radius
$1$ and balls of radius $C$ for any constant $C>0$.
\begin{rem}
Note that for graphs with bounded degree (equipped with the
counting measure), or Riemannian manifolds with bounded Ricci
curvature (equipped with the Riemannian measure), quasi-isometries
are automatically large-scale equivalences.
\end{rem}
\subsection{Volume of balls and growth function.\\}

Let $(X, d, \mu)$ be a metric measure space. The equivalence class
modulo $\approx$ of $\mu(B(x,r))$ is independent from $x$. We call
it the volume growth of $X$ and we write it $V(r)$. One can easily
show the next result (see \cite{CS}).
\begin{prop}\label{p1}
The volume growth is invariant under large-scale equivalence
(among $(DV)_{loc}$ spaces).
\end{prop}

\begin{defn}
Let $X$ be a metric measure space. We say that $X$ is doubling if
there exists a constant $C>0$ such that, $\forall x\in X$ and
$\forall r\geq 0$:
\begin{equation}\label{doubling}
\mu(B(x,2r))\leq C\mu(B(x,r)).
\end{equation}
We will call this property $(DV)$.
\end{defn}

\begin{rem}
It is easy to see that $(DV)$ is invariant under large scale
equivalence between $(DV)_{loc}$ spaces. To be more general, we
could define an asymptotic doubling condition $(DV)_{\infty}$,
restricting (\ref{doubling}) to balls of radius more than a
constant (depending on the space). Property $(DV)_{\infty}$ is
also stable under large-scale equivalence between $(DV)_{loc}$
spaces and has the advantage to focus on large scale properties
only. Actually, in every situation met in this paper, the
assumption $(DV)$ can be replaced by $(DV)_{\infty}+(DV)_{loc}$
(note that they are equivalent for graphs). Nevertheless, for the
sake of simplicity, we will leave this generalization aside.
\end{rem}

\begin{exe}
A crucial class of doubling spaces is the class of spaces with
polynomial growth: we say that a metric measure space has (strict)
polynomial growth of degree $d$ if there exists a constant $C\geq
1$ such that, $\forall x\in X$ and $\forall r\geq 1$:
$$C^{-1}r^d\leq \mu(B(x,r))\leq Cr^d.$$
Note that Gromov (see \cite{Gro}) proved that every doubling
finitely generated group is actually of polynomial growth of
integer degree (since it is virtually nilpotent).
\end{exe}

\subsection{Contents of the article.\\}

In the next section, we present a setting adapted to the study of
asymptotic isoperimetry in general metric measure spaces. The main
interest of this setting is that the ``asymptotic isoperimetric
properties" are invariant under large-scale equivalence. In
particular, it will imply that if $X$ is a $(DV)_{loc}$ and
uniformly connected (see next section) space, then the class
modulo $\approx$ of $I_h$ will not depend any more on $h$ provided
$h$ is large enough. For that reason, we will simply denote $I$
instead of $I_h$. Then, we introduce a notion of weak geodesicity
(see definition \ref{M}) which is invariant under Hausdorff
equivalence but not under quasi-isometry. We call it property (M)
since it can be formulated in terms of existence of some
``monotone" geodesic chains between any pair of points. This
property plays a crucial role when we want to obtain upper bounds
for the volume of spheres (see \cite{tes}). It will also appear as
a natural condition for some properties discussed in this paper.

Here are the two main problems concerning isoperimetry in metric
measure spaces: first, determining the asymptotic behavior of the
profile; second, finding families of subsets that optimize the
profile. The asymptotic behavior of $I$ is more or less related to
volume growth (see \cite{Coul} and \cite{P'} for the case of
finitely generated groups). In the setting of groups, the two
problems have been solved for Lie groups (and for polycyclic
groups) in \cite{PitSal} and \cite{Coul} and for a wide class of
groups constructed by wreath products in \cite{Erch}. It seems
very difficult (and probably desperate) to get general statements
for graphs with bounded degree without any regularity assumption
(like doubling property or homogeneity). On the other hand, let us
emphasize the fact that doubling condition appears as a crucial
assumption in many fields of analysis. So in this article, we will
deal essentially with doubling metric measure spaces. Without any
specific assumption on the space, balls seem to be natural
candidates for being isoperimetric subsets, especially when the
space is doubling (see Corollary \ref{co}).

One could naively think that thanks to Theorem \ref{th1},
properties like asymptotic isoperimetry of balls are stable under
large-scale equivalence. Unfortunately, it is not the case.
Indeed, Theorem \ref{th1} roughly says that if $f:X\rightarrow X'$
is a large scale equivalence between two metric measured spaces,
then for $a>0$ large enough and for every measurable subset $A$ of
$X$, the measure of the boundary of $[f(A)]_a$ is smaller than the
measure of the boundary of $A$ up to a multiplicative constant.
So, in order to apply Theorem \ref{th1} to asymptotic isoperimetry
of balls, we need the existence of some $C>0$ such that:
\begin{equation}\label{balls}
B(f(x),r-C)\subset[f(B(x,r))]_C\subset B(f(x),r+C)\quad \forall
x\in X, \forall r>0.
\end{equation}
This condition is satisfied if $f$ is a Hausdorff equivalence. But
if $f$ is a quasi-isometry, we only have:
\begin{equation}\label{balls'}
B(f(x),C^{-1}r-C)\subset[f(B(x,r))]_C\subset B(f(x),Cr+C)\quad
\forall x\in X, \forall r>0.
\end{equation}
(Note that (\ref{balls}) and (\ref{balls'}) are metric
conditions).

Let us introduce some terminology. First, let us write
$\mathcal{B}$ for the family of all balls of $X$.

\begin{defn}
Let $X$ be a metric measure space.
\begin{itemize}
\item We say that $X$ is {\bf (IB)} if balls are asymptotically
isoperimetric, i.e. if $$I^{\downarrow}_{\mathcal{B}}\preceq I.$$
Otherwise, we will say that $X$ is {\bf (NIB)}.

\item We say that $X$ is {\bf strongly-(IB)} if balls are strongly
asymptotically isoperimetric, i.e. if
$$I^{\uparrow}_{\mathcal{B}}\preceq I.$$

\item Finally, we say that a metric measure space is {\bf
stably-(IB)} (resp. stably-(NIB)) if every (M)-space (see
Definition \ref{M})  isometric at infinity to $X$ is (IB) (resp.
(NIB)). If necessary, we will restrict our study to a certain
class of metric measure spaces.
\end{itemize}
\end{defn}

\begin{defn}
We say that a space $(X,d,\mu)$ satisfies a strong (isoperimetric)
inequality---or that $X$ has a strong profile---if $I\succeq
id/\phi$ where $\phi$ is the equivalence class modulo $\approx$ of
the function:
$$t\rightarrow\inf\{r, \mu(B(x,r)\geq t\}.$$
\end{defn}

We will show that every doubling space satisfying a strong
isoperimetric inequality satisfies (IB). This actually implies
that such a space satisfies stably-(IB). In particular, any
compactly generated, locally compact group of polynomial growth
satisfies (IB). In contrast, apart from the Abelian case
\cite{tes}, it is still unknown whether such a group $G$ satisfies
strongly-(IB) or not, or, in other words, if we have
$\mu(K^{n+1}\smallsetminus K^n)\approx n^{d-1}$ where $K$ is a
compact generating set of $G$ and $\mu$ is a Haar measure on $G$.

Conversely, we will show that every strongly-(IB) doubling space
satisfies a strong isoperimetric inequality. On the other hand, we
will see that the strong isoperimetric inequality does not imply
strongly-(IB), even if the volume growth is linear ($V(r)\approx
r$).

To see that strongly-(IB) is not stable under large scale
equivalence, even among graphs with polynomial growth, we shall
construct a graph quasi-isometric to $\Z^2$ whose volume of
spheres is not dominated by $r^{\log 3/\log 2}$ (where $r$ is the
radius). Note that this can be compared with the following result
(see \cite{tes}, theorem $1$) :

\begin{thm}\cite{tes}
Let $X$ be a metric measure space with properties (M) and $(DV)$
(for instance, a graph or a complete Riemannian manifold with the
doubling property). There exists $\delta>0$ and a constant $C>0$
such that, $\forall x\in X$ and $\forall r>0$:
$$\mu\left(B(x,r+1)\smallsetminus B(x,r)\right)\leq Cr^{-\delta}\mu(B(x,r)).$$
In particular, the ratio $\mu(\partial B_{x,r}(x))/\mu(B(x,r))$
tends to $0$ uniformly in $x$ when $r$ goes to infinity.
\end{thm}

When the profile is not strong, we will see that a many situations
can happen. All the counterexamples built in the corresponding
section will be graphs of polynomial growth.

The case of a bounded profile is quite specific. Indeed, in that
case, and under some hypothesis on $X$ (including graphs and
manifolds with bounded geometry), we will prove that if
$(P_n)_{n\in \N}$ is an asymptotically isoperimetric sequence of
connected subsets of $X$, one can find a constant $C\geq 1$ and
$\forall n\in \N$, some $x_n\in X,r_n>0$ such that:
$$B(x_n,r_n)\subset P_n\subset B(x_n,Cr_n).$$
Note that here, we don't ask $X$ to be doubling.

Nevertheless, we will see that there exist graphs with polynomial
growth (with unbounded profile) such that no asymptotically
isoperimetric family has this property. In particular, those
graphs are stably-(NIB).

To be complete, we also build graphs with polynomial growth,
bounded profile and satisfying stably-(NIB).

Concerning the stability under large-scale equivalence, we will
see that even among graphs with polynomial growth, with bounded or
unbounded profile, property (IB) is not stable under large-scale
equivalence (in the case of graphs equipped with the counting
measure, a large-scale equivalence is simply a quasi-isometry).

Finally, we shall examine isoperimetric properties of connected
subsets (say that $A$ is (metrically) connected if for any
partition $A=A_1\sqcup A_2$, with $d(A_1,A_2)\geq 10$, either
$A_1$ or $A_2$ is empty). Clearly, since balls of a (M)-space are
connected, the strong isoperimetric inequality implies that
connected sets are asymptotically isoperimetric. On the other
hand, we will show that there exist graphs with polynomial growth
whose connected subsets are not asymptotically isoperimetric.

\section{Isoperimetry at infinity: a general setting}

\subsection{Isoperimetric at a given scale.\\}\label{loc}

The purpose of this section is to find some minimal conditions
under which ``isoperimetric properties at infinity" become
invariant under large-scale equivalence. In the introduction,
namely in Section \ref{bord}, we justified our definition of the
boundary by the fact that we want it to have a uniform thickness.
Nevertheless, it is not suffisant to our purpose: it is also
important for the space $X$ to look connected at scale $h$. Let
$X$ be a graph; if $h=1/2$, then every subset of $X$ has a trivial
boundary, so that all the isoperimetric properties of $X$ are
trivial.
\begin{defn}
Let $X$ be a metric space and fix $b>0$. Let us call a $b$-chain
of length $n$ from $x$ to $y$, a finite sequence
$x_0=x,\ldots,x_n=y$ such that $d(x_i,x_{i+1})\leq b$.
\end{defn}

The following definition can be used to study the isoperimetry at
a given scale, although we will only use it ``large-scale version"
in this paper.

\begin{defn}\label{localcon}

\

\noindent{\it Scaled version:} \ Let $b>0$ and $E_1\gg b$. Let us
say that $X$ is uniformly $b$-connected at scale $\leq E_1$ if
there exists a constant $E_2\geq E_1$ such that for every couple
$x,y\in X$ such that $d(x,y)\leq E_1$, there exists a $b$-chain
from $x$ to $y$ totally included in $B(x,E_2)$.

\noindent{\it Large-scale version:} \ If, for all $E_1\gg b$, $X$
is uniformly $b$-connected at scale $\leq E_1$, then we say that
$X$ is uniformly $b$-connected (or merely uniformly connected).
\end{defn}

\begin{rem}
Note that in the scaled version, the space $X$ is allowed to have
a proper nonempty subset $A$ such that $d(A,A^c)>E_1$: in this
case $X$ is not $b$-connected at all.
\end{rem}

{\bf Invariance under quasi-isometry:} Note that if $X$ is
uniformly $b$-connected at scale $\leq E_1$ and if $f:$
$X\rightarrow X'$ is a quasi-isometry of constants $C_1$ and $C_2$
smaller enough than $E$, then $X'$ is uniformly
$C_2b+C_1$-connected at scale $\leq E_1/C_2-C_1$. In particular,
if $X$ is uniformly $b$-connected, then $X'$ is uniformly
$(C_2b+C_1)$-connected.

\begin{rem}
Let us write $d_b(x,y)$ for the $b$-distance from $x$ to $y$, that
is, the minimal length of a $b$-chain between $x$ and $y$ (note
that if every couple of points of $X$ can be joined by a
$b$-chain, then $d_b$ is a pseudo-metric on $X$).

If there exists $C>0$ such that, for all $x,y\in X$, one has:
$d_b(x,y)\leq Cd(x,y)+C$, then in particular, $X$ is uniformly
$b$-connected (we can call this property: quasi-geodesic
property).
\end{rem}

\begin{exe}
A graph and a Riemannian manifold are respectively uniformly
$1$-connected and uniformly $b$-connected for all $b>0$.
\end{exe}

\begin{prop}\label{prop1}
Let $X$ be a uniformly $b$-connected space at scale $\leq E_1$.
Let $h$ be such that $h\geq 2b$.
\begin{itemize}
\item[(i)] For every subset $A$ of $X$ and every $x\in A^c$ such
that $d(x,A)< E_1$ (resp. $x\in A$ such that $d(x,A^c)< E_1$),
there exists a point $z\in \partial_h A$ at distance $\leq E_2$ of
$x$ such that:
$$B(z,b)\subset\partial_h A.$$

\item[(ii)] If, moreover, $X$ is $(DV)_{loc}$ and $h\ll E_1$, then
there exists a constant $C'\geq 1$ such that, for every subset
$A$, there exists a family $(B(y_i,b))_i$ included in $\partial_h
A$, such that, for all $i\neq j$, $d(y_i,y_j)\geq E_2$ and such
that:
$$\sum_i\mu(B(y_i,b))\leq \mu(\partial_h A)\leq C'\sum_i\mu(B(y_i,b)).$$

\item[(iii)] The $h$-boundary measure of a subset of a
$(DV)_{loc}$, uniformly $b$-connected  space does not depend on
$h$ up to a multiplicative constant, provided $E_1\gg h\geq 2b$.
\end{itemize}
\end{prop}

\bpr Let $x\in A^c$ such that $d(x,A)<E_1$ and let $y\in A$ be
such that $d(x,y)\leq E_1$. We know from the hypothesis that there
exists a finite chain $x_0=x, x_1, \ldots, x_n=y$ satisfying:
\begin{itemize}
\item $x_n\in A$,

\item $d(x,x_i)\leq E_2$ for all $i$,

\item for all $1\leq i\leq n$, $d(x_{i-1},x_i)\leq b$.

\end{itemize}
Since $x\in A^c$ and $y\in A$, there exists $j\leq n$ such that
$x_{j-1}\in A^c$ and $x_{j}\in A$. Clearly, $x_j\in A_b\cap
[A^c]_b=\partial_b A$. But since $[\partial_{b}A]_b\subset
\partial_{2b} A\subset\partial_h A$, the ball $B(x_j,b)$ is included in $\partial_h A$, which
proves the first assertion. \

Let us show the second assertion. Consider a maximal family of
disjoint balls $(B(x_i,2E_2))_{i\in I}$ with centers $x_i\in
\partial_h A$. Then $(B(x_i,5E_2))_{i\in I}$ forms a covering of $\partial_h A$.

Using the first assertion and the fact that $h\ll E_1$, one sees
that each $B(x_i,2E_2)$ contains a ball $B(y_i,b)$ included in
$\partial_h A$. It is clear that the balls $B(y_i,10E_2)$ form a
covering of $\partial_h A$ and that the balls $(B(y_i,b)$ are
disjoint. But, by property $(DV)_{loc}$, there exists $C'\geq 1$,
depending on $b$ and $E_2$, such that, for all $i\in I$:
$$\mu(B(y_i,10E_2))\leq C'\mu(B(y_i,b)).$$
We deduce: $$\sum_i\mu(B(y_i,b))\leq \mu(\partial_h A)\leq
C'\sum_i \mu(B(y_i,b))$$ which proves (ii). The assertion (iii)
now follows from (ii). \epr

\

\begin{rem}
This proposition gives conditions to study isoperimetry at scale
between $b$ and $E_1$, i.e. choosing $h$ far from those two
bounds. Thus, we will always assume that this condition holds and
we will simply write $\partial A$ instead of $\partial_h A$.
Otherwise, problems may happen. We talked about what can occur if
$h<b$ at the beginning of this section. Now, let us give an idea
of what can happen if $h>E_1$. Consider a metric measure space $X$
such that $X=\cup_{i\in I} X_i$ where the $X_i$ are subsets such
that $d(X_i,X_j)\geq E_1$ whenever $i\neq j$ and such that
$\mu(X_i)$ is finite for every $i\in I$ but not bounded. Note that
for $h<E_1$ the boundary of every $X_i$ is empty so that the
family $(X_i)_{i\in I}$ is trivially asymptotically isoperimetric.
But this can change dramatically if $h>E_1$ because the boundary
of $X_i$ can meet many $X_j$'s for $j\neq i$.
\end{rem}
\begin{rem}
If we replace uniformly $b$-connected at scale $\leq E_1$ by
uniformly $b$-connected, then the proposition gives a setting
adapted to the study of large scale isoperimetry. Namely, it says
that for a uniformly $b$-connected, $(DV)_{loc}$ space, the choice
of $h$ does not matter, provided $h\geq 2b$.
\end{rem}

\begin{cor}\label{cor2}
Let $X$ be a $(DV)_{loc}$, uniformly $b$-connected space. If
$h,h'\geq 2b$, we have:
$$I_{h}\approx I_{h'}.$$
So, from now on, we will simply call ``profile" (instead of
$h$-profile) the equivalence class modulo $\approx$ of $I_h$. Note
that the same holds for restricted profiles
$I_{h,\A}^{\downarrow}$, and $I_{h,\A}^{\uparrow}$ that we will
simply denote $I_{\A}^{\downarrow}$ and $I_{\A}^{\uparrow}$ (where
$\mathcal{A}$ is a family of subsets of $X$).
\end{cor}

The following theorem shows that a large-scale equivalence $f$
with controlled constants essentially preserves all isoperimetric
properties.

\begin{thm}\label{th1}
Let $f:$ $(X,d,\mu)\rightarrow (X',d',\mu')$ be a large-scale
equivalence (with constants $C_1$, $C_2$ and $C_3$) where $X$
(resp. $X'$) is $(DV)_{loc}$ and uniformly $b$-connected at scale
$\leq E_1$ (resp. uniformly $b'$-connected at scale $\leq E'_1$).
We suppose also that $E_1$ and $E'_1$ are far larger than $C_1$,
$C_2$, $C_2b$ and $C_2(b'+C_1)$. Then, there exists a constant
$K\geq 1$ such that, for any subset $A$ of finite measure:
$$\mu'(\partial [f(A)]_{C_1})\leq K\mu(\partial A).$$
\end{thm}

\bpr Let us start with a lemma:

\begin{lem}
Let $X$ be a $(DV)_{loc}$ space and fix some $\alpha>0$. Then
there exists a constant $c>0$ such that, for all family
$(B(x_i,\alpha))_{i\in I}$ of disjoint balls of $X$, there is a
subset $J$ of $I$ such that $\forall j\in J$, the balls
$B(x_j,2\alpha)$ are still disjoint, and such that:
$$\sum_{j\in J} \mu(B(x_j,2\alpha))\geq c\sum_{i\in I}\mu(B(x_i,\alpha)).$$
\end{lem}
\bpr Let us consider a maximal subset $J$ of $I$ such that
$(B(x_j,2\alpha))_{j\in J}$ forms a family of disjoint balls.
Then, by maximality, we get:
$$\bigcup_{i\in I} B(x_i,\alpha)\subset \bigcup_{j\in J} B(x_j,4\alpha).$$
We conclude thanks to property $(DV)_{loc}$. \epr

\

To fix ideas, take $h=2b$ and $h'=2b'$. Assertion (ii) of
Proposition \ref{prop1} implies that there exists a family of
balls $(B(y_i,b'))_i$ included in $\partial [f(A)]_{C_1}$ such
that, for all $i\neq j$, $d(y_i,y_j)\geq E_2'$ and such that:
$$\sum_i\mu(B(y_i,b'))\leq \mu(\partial_h [f(A)]_{C_1})\leq C'\sum_i\mu(B(y_i,b')).$$

By the lemma, and up to changing the constant $C'$, one can even
suppose that $d(y_i,y_j)\gg C_2E_2$ for $i\neq j$.

For all $i$, let $x_i$ be a element of $X$ such that
$d(f(x_i),y_i)\leq C_1.$ The points $x_i$ are then at distance
$\gg E_2$ to one another. Moreover, since $y_i$ is both at
distance $\leq 2b+C_1$ of $f(A)$ and of $f(A^c)$, $x_i$ is both at
distance $\ll E_1$ of $A$ and of $A^c$. So, by the assertion (i)
of the proposition, there exists a ball $B(z_i,b)$ included in
$\partial A\cap B(x_i,E_2)$. Since balls $B(x_i,E_2)$ are
disjoint, so are the $B(z_i,b)$. The theorem then follows from
property $(DV)_{loc}$ and from property of ``almost-conservation"
of the volume (property (b)) of large-scale equivalence.\epr

\

\begin{rem}
Note that in the case of graphs, the condition $h\geq 2$ can be
relaxed to $h\geq 1$ (the proposition and the theorem stay true
and their proofs are unchanged).
\end{rem}

\begin{cor}\label{cor1}
Under the hypotheses of the theorem, we have:
\begin{itemize}
\item[(i)] if the family $(A_i)_{i\in I}$ is asymptotically
isoperimetric, then so is $(f(A_i)_b)_{i\in I}$;

\item[(ii)] if $I$ and $I'$ are the profiles of $X$ and $X'$
respectively, we get: $I\approx I'$.
\end{itemize}
\end{cor}

The corollary results immediately from the theorem and the
following proposition. \epr

\begin{prop}\label{propo2}
Let $f$ be a large-scale equivalence between two $(DV)_{loc}$
spaces $X$ and $X'$. Then for every subset $A$ of $X$, there
exists $C\geq 1$ such that:
$$\mu(A)\leq C\mu'([f(A)]_{C_1}).$$
\end{prop}
\bpr Consider a maximal family of disjoint balls
$(B(y_i,C_1))_{i\in I}$ whose centers belong to $f(A)$. These
balls are clearly included in $[f(A)]_{C_1}$. By property
$(DV)_{loc}$, the total volume of these balls, and therefore
$\mu'([f(A)]_{C_1})$, are comparable to the sum of the volumes of
balls $B(x_i,3C_1)_{i\in I}$ that form a covering of
$[f(A)]_{C_1}$. The preimages of these balls thus cover $A$. But,
for each $i$, $f^{-1}(B(y_i,3C_1))$ is contained in a ball of
radius $3C_1C_2+C_2$ and of center $x_i$ where $x_i\in
f^{-1}(\{y_i\})$. By property $(DV)_{loc}$ and property of
almost-conservation of the measure of small balls (property (b))
of $f$, the measure of this ball is comparable to that of
$B(y_i,3C_1)$. So we are done. \epr

\

Finally, let us mention that if we suppose that $X$ and $X'$ are
uniformly connected and satisfy the $(DV)_{loc}$ condition, then
Theorem \ref{th1} and its corollary hold for any large-scale
equivalence $f$.

\subsection{Property (M): monotone geodesicity.\\}

Let us introduce a natural (but quite strong) property of
geodesicity.
\begin{defn}\label{M}
Let us say that $(X,d)$ has property (M) if there exists $C\geq 1$
such that, $\forall x\in X$, $\forall r>0$ and $\forall y\in
B(x,r+1)$, we have $d(y,B(x,r))\leq C.$
\end{defn}

\begin{rem}
Let $(X,d)$ be a (M) metric space. Then $X$ has ``monotone
geodesics" (this is why we call this property (M)): i.e. there
exists $C\geq 1$ such that, for all $x,y\in X$, there exists a
finite chain $x_0=x, x_1,\ldots,x_n=y$ such that $\forall 0\leq
i<n$,
$$d(x_i,x_{i+1})\leq C;$$
and
$$d(x_{i},x)\leq d(x_{i+1},x)-1.$$
Consequently, $\forall r,k>0$, $\forall y\in B(x,r+k)$, we have:
$$d(y,B(x,r))\leq Ck.$$
These two properties are actually trivially equivalent to property
(M).
\end{rem}

Recall (see \cite{Grom}, p 2) that two metric spaces $X$ and $Y$
are said Hausdorff equivalent:
$$X\sim_{Hau}Y$$
if there exists a (larger) metric space $Z$ such that $X$ and $Y$
are contained in $Z$ and such that
$$\sup_{x\in X}d(x,Y)<\infty$$
and
$$\sup_{y\in Y}d(y,X)<\infty.$$
\begin{rem}
It is easy to see that property (M) is invariant under Hausdorff
equivalence. But on the other hand, property (M) is unstable under
quasi-isometry. To construct a counterexample, one can
quasi-isometric embed $\R_+$ into $\R^2$ such that the image,
equipped with the induced metric does not have property (M):
consider a curve starting from $0$ and containing for every $k\in
\N$ a half-circle of radius $2^k$. So it is strictly stronger than
quasi-geodesic property (\cite{Grom}, p 7), which is invariant
under quasi-isometry: $X$ is quasi-geodesic if there exist two
constants $d>0$ and $\lambda>0$ such that for all $(x,y)\in X^2$
there exists a finite chain of points of $X$:
$$x=x_0,\ldots,x_n=y,$$
such that
$$d(x_{i-1},x_{i})\leq d,\quad i=1\ldots n,$$
and
$$\sum_{i=1}^n d(x_{i-1},x_i)\leq \lambda d(x,y).$$
\end{rem}

\begin{exe}
A geodesic space has property (M), in particular graphs and
complete Riemannian manifolds have property (M). A discretisation
(i.e. a discrete net) of a Riemannian manifold $X$ has property
(M) for the induced distance.
\end{exe}

\begin{rem}\label{re1}
Note that in general, if $X$ is a metric measure space, we have:
$$\partial_{1/2}B(x,r+1/2)\subset B(r+1)\setminus
B(x,r).$$ Moreover, if $X$ has property (M), then, we have:
$$B(x,r+1)\setminus B(x,r)\subset \partial_C B(x,r+1).$$
\end{rem}
Note that this is not true in general, even for quasi-geodesic
spaces.

\section{Link between isoperimetry of balls and strong isoperimetric inequality}

\subsection{Strong isoperimetric inequality implies (IB)}

\

To fix ideas, spaces we will consider from now on will be
$(DV)_{loc}$ and uniformly $1$-connected. Let us write: $\partial
A=\partial_{2}A$ for any subset $A$ of a metric space $X$.

\

Let $X$ be a metric measure space. Let $V$ be a nondecreasing
function belonging to the volume growth class (for instance:
$V(r)=\mu(B(x,r))$ for a $x\in X$). Write
$\phi(t)=\inf\{r,V(r)\geq t\}$ for the ``right inverse" function
of $V$. Remark that if $f$ and $g$ are nondecreasing functions
$\R_+\rightarrow\R_+$, then $f\approx g$ if and only if their
right inverses are equivalent. In particular, the equivalence
class of $\phi$ is invariant under large-scale equivalence.

\begin{defn}
Let us call a strong isoperimetric inequality the following kind
of isoperimetric inequality:
$$\forall A\subset X, \quad |\partial A|\geq  C^{-1}|A|/\phi(C|A|).$$
Remark that this is equivalent to:
$$I\succeq id/\phi,$$
Therefore, if $X$ satisfies a strong isoperimetric inequality, we
will say that it has a strong profile.
\end{defn}

\begin{exe}
If $X$ has polynomial growth of degree $d$, we have:
$\phi(t)\approx t^{1/d}$. So $X$ has a strong profile if and only
if:
$$I\succeq (id)^{\frac{d-1}{d}}.$$
\end{exe}
Write, for all $x\in X$ and for all $0<r<r'$:
$$C_{r,r'}(x)=B(x,r')\setminus B(x,r).$$

\begin{prop}\label{p}
Let $X$ be a doubling space (here, no other hypothesis is
required). There exists a constant $C\geq 1$ such that:
$$\forall x\in X,\forall r\geq 1,\quad \inf_{r\leq r'\leq 2r}\mu(C_{r'-1,r'})\leq
C\mu(B(x,r))/r.$$
\end{prop}

\bpr Clearly, it suffices to prove the proposition when $r=n$ is a
positive integer. First, note that
$$\cup_{k=n}^{2n}(B(x,k)\smallsetminus B(x,k-1))\subset B(x,2n).$$
So, we have
$$\mu(B(x,2n))\geq n\inf_{n \leq k\leq 2n}\mu(B(x,k)\setminus
B(x,k-1)).$$ We conclude by Doubling property. \epr

\

\begin{cor}\label{co}
Let $X$ be a uniformly connected doubling space. Then we have:
$$I^{\downarrow}_{\mathcal{B}}\preceq id/\phi.$$ Namely, there exists a constant $C\geq 1$ such that:
$$\forall x\in X,\forall r>0,\quad \inf_{r'\geq r}\mu\left(\partial B(x,r')\right)\leq
C\mu(B(x,r))/r.$$
\end{cor}
{\it Proof.} This follows from Remark \ref{re1}. \epr

\

\begin{cor}\label{p3}
Let $X$ be a uniformly connected doubling space satisfying a
strong isoperimetric inequality. Then, $X$ is stably-(IB).
\end{cor}
\bpr It follows from Corollary \ref{co} and from Corollary
\ref{cor1}. \epr

\begin{rem}
Varopoulos \cite{Varo} showed that the strong isoperimetric
inequality is satisfied by any group of polynomial growth. Coulhon
and Saloff-Coste \cite{Coul} then proved it for any unimodular
compactly generated locally compact group with a simple and
elegant demonstration. We have the following corollary.
\end{rem}

\begin{cor}
A Cayley graph of a group of polynomial growth is stably-(IB).
\end{cor}

\subsection{The strong isoperimetric inequality does not imply
strongly-(IB)}\label{s2}

\

Note that this will result from the example shown in section \ref{s1}.
Let us present here a counterexample with linear growth.

For every integer $n$, we consider the following finite rooted
tree $G_n$: first take the standard binary tree of depth $n$. Then
stretch it as follows: replace each edge connecting a $k-1$'th
generation vertex to a $k$'th generation vertex by a (graph)
interval of length $2^{2^{n-k}}$. Then consider the graph $G_n'$
obtained by taking two copies of $G_n$ and identifying the
vertices of last generation of the first copy with those of the
second copy. Write $r_n$ and $r'_n$ for the two vertices of $G_n'$
corresponding to the respective roots of the two copies of $G_n$.
Finally, glue ``linearly" the $G_n'$ together identifying $r_n'$
with $r_{n+1}$, for all $n$: it defines a graph $X$.

\

Let us show that $X$ has linear growth (i.e. polynomial growth of
degree $1$). Thus $I\approx 1$, and since the boundary volume of
balls is clearly not bounded, we do not have
$I_B^{\uparrow}\preceq I$. In particular, $X$ is not
strongly-(IB).

Since $X$ is infinite, it is enough to show that there exists a
constant $C>0$ such that
\begin{equation}\label{eq'}
|B(x,r)|\leq Cr
\end{equation}
for every vertex $x$ of $X$. But it is clear that among the balls
of radius $r$, those which are centered in points of $n$'th
generation of a $G_n$ for $n$ large enough are of maximal volume.
Let us take such an $x$. Remark that for
$\sum_{j=0}^{k}2^{2^j}\leq r\leq \sum_{j=0}^{k+1}2^{2^{j}}$, we
have:
$$|B(x,r)|\leq 2\mid B(x,\sum_{j=0}^{k}2^{2^j})\mid+2\left(r-\sum_{j=0}^{k}2^{2^j}\right)$$
So it is enough to show (\ref{eq'}) for $r=\sum_{j=0}^{k}2^{2^j}$.
We have:
$$\mu(B(x,\sum_{j=0}^{k}2^{2^j}))=\sum_{j=0}^{k}2.2^j.2^{2^{k-j}}\leq 4.2^{2^k}.$$
Which proves ($\ref{eq'}$) with $C=8$. \epr

\

\begin{rem}\label{r''}
This example and that of section \ref{s1} show in particular that
the strong isoperimetric inequality does not imply (even in linear
growth case) strongly-(IB).
\end{rem}

\subsection{Instability of strongly-(IB) under quasi-isometry}\label{s1}

\begin{thm}
We can find a graph, quasi-isometric to $\Z^2$ (resp. a Riemannian
manifold $M$ bi-Lipschitz equivalent to $\R^2$) whose volume of
spheres is not dominated by $r^{\log 3/\log 2}$ (where $r$ is the
radius).
\end{thm}
\begin{rem}
The restriction to dimension $2$ is not essential, but was made to
simplify the exposition (actually, we merely need the dimension to
be greater or equal to $2$).
\end{rem}

\bpr The general idea of the construction is to get a sequence of
spheres which look like finitely iterated Von Koch curves. First,
we will build a graph with weighted edges. Actually, this graph
will be simply the standard Cayley graph of $\Z^2$, and the edges
will have lengths equal to $1$ except for some selected edges
which will have length equal to a small, but fixed positive
number.

\noindent{\it First step of the construction:} Let us define a
sequence $(A_k)$ of disjoint subtrees of $\Z^2$ (which is
identified to its usual Cayley graph). Let $(e_1,e_2)$ be the
canonical basis of $\Z^2$ and denote $S=\{\pm e_1,\pm e_2\}$. For
every $k\geq 1$, let $a_k=(2^{2k},0)$ be the root of the tree
$A_k$ and define $A_k$ by:
\begin{equation}\label{tree}
x\in A_k\Leftrightarrow
x=a_k+2^{k}\eps_0(x)+2^{k-1}\eps_1(x)+\ldots+2^{k-i(x)}\eps_{i(x)}(x)+r(x)\eps_{i(x)+1}(x)
\end{equation}
where:

\noindent--- $0\leq i(x)\leq k-1,$

\noindent--- $\eps_j(x)$ belongs to $S$ for every $0\leq j\leq
i(x)+1$ and is such that $\eps_{j+1}(x)\neq -\eps_j(x)$ (for
$j\leq i(x)$),

\noindent--- $r(x)\leq 2^{k-i(x)-1}-1$.

It is easy to see that $A_k$ is a subtree of $\Z^2$ and that the
above decomposition of $x$ is unique. In particular, we can
consider its intrinsic graph metric $d_{A_k}$: let $S_k$ be the
sphere of center $a_k$ and of radius $2^{k+1}-1$ for this metric.
Clearly, $|S_k|\geq 3^{k-1}.$

\noindent{\it Second step of the construction:} We define a graph
$Y$ with weighted edges as follows: $Y$ is the usual Cayley graph
of $\Z^2$; all edges of $Y$ have length $1$ but those belonging to
$A=\cup_k A_k$ which have length equal to $1/100$. The measure on
$Y$ is the countable measure and the distance between two vertices
$v$ and $w$ is the minimal length of a chain joining $v$ to $w$,
the length of a chain being the sum of the weights of its edges.
Clearly, as a metric measure space, $Y$ is large-scale equivalent
to $\Z^2$.

For every $k\geq 2$, consider the sphere
$S(a_k,r_k)=B(a_k,r_k+1)\smallsetminus B(a_k,r_k)$ of $Y$, where
$r_k=(2^{k+1}-1)/100$.

\begin{clai}
We have: $S_k\subset S(a_k,r_k)$, so that:
$$\mu(S(a_k,r_k))\geq 3^{k-1}\geq r_k^{\log 3/\log 2}.$$
\end{clai}
\bpr Note that the claim looks almost obvious on a drawing.
Nevertheless, for the sake of completeness, we give a
combinatorial proof. Let us show that a geodesic chain in the tree
$A_k$ is also a minimizing geodesic chain in $Y$. Applying this to
a geodesic chain between $a_k$ and any element of $S_k$ (which is
of length $r_k$ in $Y$), we have that $S_k\subset S(a_k,r_k)$, so
we are done.

So let $x$ be a vertex of $A_k$. By (\ref{tree}), we have:
$$x=a_k+2^{k}\eps_0(x)+2^{k-1}\eps_1(x)+\ldots+2^{k-i(x)}\eps_{i(x)}(x)+r(x)\eps_{i(x)+1}(x)$$
Let us show by recurrence on $d_Y(a_k,x)$ (which takes discrete
values) that
$$d_Y(a_k,x)=d_{A_k}(a_k,x)/100=(2^{k}+\ldots+2^{k-i(x)}+r(x))/100=\frac{2^{k+1}(1-2^{-i(x)-1}+r(x))}{100}$$
If $x=a_k$, there is nothing to prove. Consider
$c=(c(0)=x,c(2),\ldots,c(m)=a_k)$ a minimal geodesic chain in $Y$
between $a_k$ and $x$. Clearly, it suffices to prove that
$c\subset A_k$. Suppose the contrary. Let $t$ be the largest
positive integer such that $c(t)$ belongs to $A_k$ and $c(t+1)$
does not. Let $l$ be the smallest positive integer such that
$c(t+l)\in A_k$, so that $(c(t+1),\ldots, c(t+l-1))$ is entirely
outside of $A_k$. By recurrence, the chain $(c(t+l),\ldots,c(m))$
is in $A_k$. Thus we have:
$$d_Y(x,a_k)=d_{A_k}(x,c(t))/10+|c(t)-c(t+l)|_{\Z^2}+d_{A_k}(c(t+l),a_k)/100.$$
Since $c$ is a minimal chain, we also have:
$$d_Y(c(t),a_k)=|c(t)-c(t+l)|_{\Z^2}+d_{A_k}(c(t+l),a_k)/100.$$
The following lemma applied to $u=c(t)$ and $v=c(t+l)$ implies
that $t=t+l$ which is absurd since it means that $c$ is included
in $A_k$. \epr

\begin{lem}
let $u$ and $v$ be in $A_k$. We have:
$$|u-v|_{\Z^2}\geq (d_{A_k}(u,a_k)-d_{A_k}(v,a_k))/50.$$
\end{lem}

\bpr We can of course assume that $d_{A_k}(u,a_k)\geq
d_{A_k}(v,a_k)$. Let $u=u_1+u_2$ and $v=v_1+v_2$ with
$$u_1=2^{k}\eps_0(u)+\ldots 2^{k-i(v)}\eps_{i(v)}(u)$$
and
$$v_1=2^{k}\eps_0(v)+\ldots 2^{k-i(v)}\eps_{i(v)}(v).$$
Note that by construction,
$$d_{A_k}(u_1,a_k)=d_{A_k}(v_1,a_k)$$
and since $A_k$ is a tree,
\begin{equation}\label{diff}
d_{A_k}(u,a_k)-d_{A_k}(v,a_k)=d_{A_k}(u_2,a_k)-d_{A_k}(v_2,a_k)\leq
2^{k-i(v)+2}.
\end{equation}
On the other hand, we have:
$$|u-v|_{\Z^2}\geq ||u_1-v|_{\Z^2}-|u_2-v|_{\Z^2}|$$
First, assume that $u_1\neq v_1$. Then, by (\ref{tree}), the
projection of $u_1-v_1$ along $e_1$ or $e_2$ is not zero and
belongs to $2^{k-i(v)}\N$. Moreover, using the fact that
$\eps_{j+1}(u)\neq -\eps_{j}(u))$ for every $j$, the same
projection of $u_2-v_2$ is (in $\Z^2$-norm) less than
$$2.(2^{k-i(v)-2}+2^{k-i(v)-4}+\ldots=2^{k-i(v)-1}(1+1/4+1/4^2+\ldots)\leq 2/3.2^{k-i(v)}$$
Thus,
$$|u-v|_{\Z^2}\geq 2^{k-i(v)}/3.$$
So we are done.

Now, assume that $u_1=v_1$. If $i(u)=i(v)$ or if $i(u)\leq i(v)+1$
and $\eps_{i(v)+1}(u)=\pm\eps_{i(v)+1}(v)$, then we have
trivially:
$$|u-v|_{\Z^2}=(d_{A_k}(u,a_k)-d_{A_k}(v,a_k)).$$
Otherwise, we have:
$$u-v=u_2-v_2=(2^{k-i(v)-1}-r(v))\eps_{i(v)+1}(u)+2^{k-i(v)-2}\eps_{i(v)+2}(u)+\ldots +r(u)\eps_{i(u)+1}.$$
So, projecting this in the direction of $\eps_{i(v)+2}(u)$, and
since $\eps_{i(v)+3}(u)\neq -\eps_{i(v)+2}(u)$, we obtain:
$$|u-v|_{\Z^2}=|u_2-v_2|_{\Z^2}\geq 2^{k-i(v)-2}-(2^{k-i(v)-4}+\ldots+2^{k-i(u)}+r(u))\geq 2^{k-i(v)-2}-2^{k-i(v)-3}=2^{k-i(v)-3}.$$
Together with \ref{diff}, we get:
$$|u-v|_{\Z^2}\geq 32(d_{A_k}(u,a_k)-d_{A_k}(v,a_k))$$
which proves the lemma. \epr

\

Clearly, $Y$ is quasi-isometric to $\Z^2$. It is not difficult
(and left to the reader) to see that we can adapt the construction
to obtain a graph.

Now, let us explain briefly how we can adapt the construction to
obtain a Riemannian manifold bi-Lipschitz equivalent to $\R^2$.
First, we embed $\Z^2$ into $\R^2$ in the standard way, so that
$A_k$ is now a subtree of $\R^2$. Let $\tilde{A}$ be the
$1/100$-neighborhood of $A$ in $\R^2$. Let $f$ be a nonnegative
function defined on $\R^2$ such that: $1-f$ is supported by
$\tilde{A}$, $f\geq a$ and $f(x)=a$ for all $x\in A$. Finally,
define a new metric on $\R^2$ multiplying the Euclidian one by
$f$. \epr

\subsection{Strongly-(IB) implies the strong isoperimetric inequality.\\}

The converse to Proposition \ref{p3} is clearly false (see the
examples of the next section). However, one has:

\begin{prop}\label{p'}
Let $X$ be a doubling (M)-space. Suppose moreover that there
exists $x\in X$ such that the family of balls of center $x$ is
strongly asymptotically isoperimetric. Then we have
$$I^{\downarrow}_{\mathcal{B}}\succeq id/\phi.$$ In particular, $X$ satisfies a strong isoperimetric inequality.
\end{prop}

{\it Proof.} Since $(B(x,r))_r$ forms an asymptotically
isoperimetric family, it is enough to show that there exists $c>0$
such that:
$$\mu(\partial B(x,r))\geq c\frac{\mu(B(x,r))}{r}.$$

But, let us recall that property (M) implies that there exists
$C>0$ such that, for all $r>0$:
$$\mu(C_{r,r+1}(x))\leq C\mu(\partial_1 B(x,r)).$$
Since $(B(x,r))_r$ forms an asymptotically isoperimetric family,
there exists $C'\geq 1$, such that, for all $r'<r$:
$$\mu(\partial B(x,r'))\leq C'\mu(\partial B(x,r)).$$
Using these two remarks, we get:
$$\mu(B(x,r))\leq CC'r\mu(\partial B(x,r)).$$
So we are done. \epr

\section{What can happen if the profile is not strong}

All the metric measure spaces built in this section will be graphs
with polynomial growth. For simplicity, we write $|A|$ for the
cardinal of a finite subset $A$ of a graph.

\subsection{Bounded profile: connected isoperimetric sets are ``controlled" by balls.\\}

We will say that a subset $A$ of a metric space is metrically
connected (we will merely say ``connected" from now on) if there
does not exist any nontrivial partition of $A=A_1\sqcup A_2$ with
$d(A_1,A_2)\geq 10$.

Let $X$ be a uniformly $1/2$-connected space, with bounded
profile, and such that balls measure of radius $1/2$ is more than
a constant $a>0$. Actually, we can ignore nonconnected sets.
Indeed if $(A_n)$ is an isoperimetric family, then the $A_n$ have
a bounded number of connected components: otherwise, by
Proposition \ref{p1}, the boundary of $A_n$ would not be bounded
(because the distinct connected components have disjoint
$1$-boundaries each one containing a ball of radius $1/2$). It is
enough to replace $A_n$ by its connected component of maximal
volume.

\begin{clai}
Let $(X,d,\mu)$ be a $(DV)_{loc}$, uniformly $1/2$-connected space
such that the measure of balls of radius $1/2$ is more than $a>0$
and whose profile $I$ is bounded. Then, if $(A_n)$ is an
isoperimetric sequence of connected subsets of $X$, there exist a
constant $C>0$, some $x_n\in X$ and some $r_n>0$ such that:
$$\forall n,\quad B(x_n,r_n)\subset A_n\subset B(x_n,Cr_n).$$
\end{clai}
\bpr To fix ideas, let us assume that $\partial A=\partial_1 A$
(for all $A\subset X$). Let $y_n$ be a point of $A_n$ and write
$d_n=\sup_{y\in
\partial A_n}d(y_n,y)$. Let $r\leq d_n$ be such that $C_{r,r+1}(y_n)$
intersects nontrivially $\partial A_n$ (recall that
$C_{r,r'}(x)=B(x,r')\setminus B(x,r)$). Then, by Proposition
\ref{prop1}, there exists a constant $C\geq 1$ such that
$C_{r-C,r+C}(y_n)\cap\partial A_n$ contains a ball of radius $1/2$
and therefore has measure $\geq a$. Consequently, if
$\delta_n=\sup\{r'-r; C_{r,r'}(y_n) \cap \partial
A_n=\emptyset\}$, then:
\begin{equation}\label{equa1}
\mu(\partial A_n)\geq \frac{d_n}{2C\delta_n}a.
\end{equation}

Since the boundary of $A_n$ has bounded measure, there exists a
constant $c>0$ and, for all $n$, two positive reals $r'_n$ and
$r"_n$ such that $r"_n-r'_n\geq cd_n$ and
$C_{r'_n,r"_n}\cap\partial A=\emptyset.$

Write $s_n=(r'_n+r"_n)/2.$ Since $A_n$ is connected,
$C_{s_n-10,s_n+10}(x)\cap A_n$ is nonempty. But then, if $x_n\in
C_{s_n-10,s_n+10}(x)\cap A_n$, we get:
$$B\left(x_n,\frac{r"_n-r'_n}{2}-10\right)\subset A_n.$$
On the other hand:
\begin{equation}\label{equa2}
A_n\subset B(x_n,2d_n).
\end{equation}
Write: $r_n=cd_n/2-10$. The proposition follows from (\ref{equa1})
and from (\ref{equa2}). \epr

\subsection{Stably-(NIB) graphs with unbounded profile
and where isoperimetric families can never be ``controlled" by
families of balls}

\begin{thm}
For every integer $d\geq 2$, there exists a graph $X$ of
polynomial growth of degree $d$, with unbounded profile,
satisfying stably-(NIB) and such that, for all isoperimetric
sequences $(A_n)$, it is impossible to find sequences of balls
$B_n=B(x_n,r_n)$ and $B'_n=B(x'_n,r'_n)$ of comparable radii (i.e.
such that $r'_n/r_n$ is bounded) such that:
$$B_n\subset A_n\subset B'_n,\quad \forall n.$$
\end{thm}

Consider the graph $X$ obtained from $\Z^d$ taking off some edges.
Consider, in the axis $\Z.e_1$, the intervals $(S_n)$ of length
$[\sqrt{n}]$ and at distance $2^n$ from one another. Consider the
sequence $(A_n)$ of full parallelepiped defined by the equations:
$x_1\in I_n$ and $|x_i|\leq n/2$ for $i\geq 2$.

Then consider a partition of the boundary (in $\Z^d$) of $A_n$ in
$(d-1)$-dimensional cubes $a_n^k$ whose edges have length
approximatively $\sqrt{n}$. Remove all the edges that connect
$A_n$ to its complementary but those connected to the ``center" of
$a_n^k$ (here, the center of $a_n^k$ is a point of $\Z^d$ we
choose at distance $\leq 2$ from the ``true center" in $\R^n$ of
the convex hull of $a_n^k$). We thus obtain a connected graph $X$.
Note that the $A_n$ are such that
$$|A_n|\approx n^{d-1}\sqrt{n}$$
and:
$$|\partial_{X} A_n|\approx \frac{|\partial_{\Z^d}A_n|}{|a_{n}^0|}\approx n^{d-1}/(\sqrt{n})^{d-1}=(\sqrt{n})^{d-1}.$$

Write $A$ for the union of $A_i$ and $A^c$ for its complementary
in $X$.

\begin{clai}
The growth in $X$ is polynomial of degree $d$.
\end{clai}
\bpr It will follow from the strong profile of balls. \epr

\begin{clai}
The profile of $X$ is not strong.
\end{clai}
\bpr Let us consider the $A_n$. If the profile was strong, the
sequence  $u_n=\frac{|A_n|}{|\partial A_n|^{\frac{d}{d-1}}}$ would
be bounded. But there exists a constant $c>0$ such that:
$$u_n\geq cn^{d-1}\sqrt{n}/(\sqrt{n})^d= cn^{\frac{d-1}{2}}\rightarrow\infty.$$
\epr

\

\begin{clai}\label{f1}
Let $R$ be a unbounded subset of $\R_+$ and let $(P_r)_{r\in R}$
be a family of subsets such that there exist two constants $C\geq
1$ and $a>0$ such that:
$$\forall r>0,\exists x_r\in X, \quad B(x_r,r/C)\subset [P_r]_a\subset B(x_r,Cr).$$
Then there exists a constant $c'$ such that:
$$\forall r>0, \quad \mu(\partial P_r)\geq c'\mu(P_r)^{\frac{d-1}{d}}.$$
\end{clai}
The following lemma and its proof will be useful in all examples
that we will expose in the following sections. Write $A^c$ for the
complementary of $A$ (in $X$ or, which is actually the same in
$\Z^d$).

\begin{lem}
The profile of $A^c$ (or of $A'^c$) is strong. That means:
$I(t)\approx t^{\frac{d-1}{d}}$.
\end{lem}\label{l1}
{\bf Proof of the lemma.}

First of all, it is enough to consider only connected subsets $P$
of $A^c$. Indeed, if $P$ has many connected components $P_1\ldots
P_k$, then, by subadditivity of the function $\phi:t\rightarrow
t^{\frac{d-1}{d}}$, if the $P_i$ verify: $|\partial P_i|\geq
c\phi(|P_i|)$, then so do $P$.

Note that $A^c$ embeds into $X$ and into $\Z^d$. The idea consists
in comparing the profile of $A^c$ to that of $\Z^d$. First of all,
let us assume that a connected subset $P$ of $A^c$---seen in
$X$---intersects the boundary of many $A_n$. Then, as $|A_n|$ is
negligible compared to the distance between the $A_n$ when $n$
goes to infinity, the set of points of $\partial_{\Z^d}P$ at
distance $1$ of $A$ has negligible volume compared to $|\partial
P|$. Thus, if $|P|$ et $n$ are large enough, we get:
$$|\partial_{A^c} P|\geq \frac{1}{2}|\partial_{\Z^d} P|.$$

So it is enough to consider subsets meeting only one $A_n$. But
the complementary of a convex polyhedra of $\Z^d$ has trivially
the same profile (up to a constant) as $\Z^d$. So we are done.
\epr

\

 {\bf Proof of the claim \ref{f1}.}
Let $(P_r)$ be a family of subsets of $X$ satisfying the condition
of the proposition. We have to show that $\forall r$, $|\partial
P_r|\geq c'|P_r|^{\frac{d-1}{d}}$. If $P\subset A^c$, the claim is
a direct consequence of the lemma.

Suppose that $P$ meets some $A_n$ and that $r\geq 100C\sqrt{n}$.
Then we have already seen (in the proof of Lemma \ref{l1}) that if
many $A_n$ intersect $P_r$, the cardinal of the intersection of
this $P_r$ with $A$ are negligible compared to its boundary
provided $n$ and $|P_r|$ are large enough. We can thus suppose
that $P_r$ meets only one $A_n$. Furthermore, since $r\geq
100\sqrt(n)$, there is some $x'$ in $B(x_r,r/C)$ such that
$$B(x',r/10C)\in B(x_r,r/C)\cap A^c.$$
Then, observe that since
$B(x',r/10C)\subset [P_r]_a$, there is a constant $c>0$ such that:
\begin{equation}\label{equa3}
|P_r\cap B(x',r/C)|\geq c|B(x',r/C)|.
\end{equation}
It follows that the intersection of $P_r$ with $A^c$ has volume
$\geq c'|P_r|$ where $\bar{c}$ is a constant depending only on $C$
and $a$. So by Lemma \ref{l1}, we have:
$$|\partial_X P_r|\geq |\partial_{A^c}(P_r\cap A^c)|\geq c|P_r|^{\frac{d-1}{d}}.$$

We then have to study the case $r\leq 100C \sqrt{n}$. We can
assume that $x_r\in A_n$ (otherwise, we conclude with Lemma
\ref{l1}). Let $\pi$ be the orthogonal projection on the
hyperplane $x_2=0$. Then for $n$ large enough, $Cr$ is smaller
than $n/2$. Consequently, since $P_r\in B(x_r,Cr)$, every point of
$\pi(P_r)$ has at least one antecedent in the boundary of $P_r$.
So, we have:
$$|\partial_X P_r|\geq |\pi(P_r)|.$$
Moreover, note that $\pi(B(x_r,r/C))=B(\pi(x_r),r/C)$ (note that
this ball lies in $\Z^{d-1}$). On the other hand, since the
projection is $1$-Lipschitz, we get:
$$\pi([P_r]_a)\subset [\pi(P_r)]_a,$$
so
$$B(\pi(x_r),r/C)\subset [\pi(P_r)]_a.$$
Similarly to (\ref{equa3}), we have:
$$|\pi(P_r)\cap B(\pi(x_r),r/C|\geq c|B(\pi(x_r),r/C)|$$
So, finally, we have:
$$|\partial_X P_r|\geq c'r^{d-1}$$
so we are done. \epr

\

\begin{cor}
In every space isometric at infinity to $X$, the volume of spheres
$\approx r^{d-1}$. In particular, they are not asymptotically
isoperimetric.
\end{cor}

\textit{Proof of the corollary.} Let $f$: $X'\rightarrow X$ a
large-scale equivalence between two metric measure spaces $X'$ and
$X$ and take $y\in X'$. It comes:
$$B\left(f(y),\frac{r}{C_2}-C_1\right)\subset B([f(B(y,r))]_{C_1})\subset
B(f(y),C_2r+C_1).$$ The corollary follows from Claim \ref{f1} and
from Theorem \ref{th1}. \epr

\subsection{Graphs stably-(NIB) with bounded profile}

\begin{thm}
For any integer $d\geq 2$, one can find graph of polynomial growth
of degree $d$, with bounded profile, and which is stably-(NIB).
\end{thm}

The construction follows the same lines as in the previous
section. Consider in $\Z^d$, a sequence $(C_n)$ of subsets defined
by:
$$C_n=B(x_n,n)\cup B(x'_n,n)$$
where $x_n=(2^{n+1},n-\log_n,0,\ldots,0)$ and
$x'_n=(2^{n+1},\log_n-n,0,\ldots,0)$.

We disconnect $C_n$ from the rest everywhere but in the axis
$\Z.e_1$. Let $Y$ be the corresponding graph. $C_n$ looks like a
ball (of $\Z^d$) ``constricted" at the equator. Indeed, every
point of $C_n$ belonging to the hyperplane $\{x_2=0\}$ is at
distance at most $\log n$ from the boundary (in $Y$) of $C_n$.
This is the property that will prevent $C_n$ from being
``deformed" into a ball. Write $C=\cup_n C_n$.

\begin{lem}\label{l2}
The graph $C^c$ has a strong profile.
\end{lem}
\bpr The demonstration is essentially the same as for Lemma
\ref{l1}. \epr

\begin{clai}
The growth in the graph $X$ is polynomial of degree $d$.
\end{clai}
\bpr We have to show that there exists a constant $c>0$ such that,
$\forall x,r,$  $|B(x,r)|\geq cr^d$ (the converse inequality
following from the fact that $X$ embeds in $\Z^d$). Thanks to
Lemma \ref{l2}, we can suppose that $B$ is included in a $C_{n_0}$
so that its radius is $\leq n_0$.

The conclusion follows then from the next trivial fact: in $\Z^d$,
if $r\leq n_0$, the volume of the intersection of a ball of radius
$n_0$ with a ball of radius $r\leq n_0$ and of center belonging to
the first ball is $\geq 2^{-d}|B(x,r)|\geq 2^{-10d}r^d$. Indeed,
the worst case is when $x$ is in a ``corner" of the ball. So we
are done. \epr

\

\begin{clai}\label{f2}
If $Y'$ is a (M)-space which is isometric at the infinity to $Y$,
then its balls are not asymptotically isoperimetric.
\end{clai}
\bpr The demonstration results from the following lemma and
Proposition~\ref{p1}.

\

\begin{lem}
Let $\mathbf{P}$ be an asymptotically isoperimetric family of
connected subsets of $X$. Then there exists a constant $C\geq 1$
such that, for all $P\in\mathbf{P}$ of measure $>C$, there exists
$n$ such that $|P\bigtriangleup C_n|\leq C$.
\end{lem}
\bpr Since the profile of $C^c$ is strong, it is clear that for
$|P|$ large enough, $P\cap C^c$ must be bounded. We then have to
show that if $(P_n)$ is a sequence of subsets such that for all
$n$, $P_n\subset C_n$ and such that $|P_n|$ and $|C_n\setminus
P_n|$ tends to infinity, then $|\partial P_n|$ also tends to
infinity. Suppose, for instance that $|P_n|\leq |C_n\setminus
P_n|$. But Theorem \ref{th1} makes clear that this problem in
$\Z^d$ is equivalent to the similar problem in $\R^d$: that is,
replacing $C_n$ with its convex hull $\tilde{C_n}$ in $\R^d$.
Since the $\tilde{C}_n$ are homothetic copies of $\tilde{C_1}$, by
homogeneity, we only have to show that the profile $I(t)$ of
$\tilde{C_1}$ is $\geq ct^{\frac{d-1}{d}}$ for
$0<t<|\tilde{C_1}|/2$, which is a known fact (see \cite{Ros}).
\epr

\

Let us finish the demonstration of Claim \ref{f2}. We now have to
show that the sets $C_n$ cannot be---up to a set of bounded
measure---inverse images of balls by some large-scale equivalence.
So let $(X',d,\mu)$ be a (M)-space and let $f:$ $X\rightarrow X'$
be a large-scale equivalence.

Let us consider two points of $C_n$ of respectively maximum and
minimum $x_2$. The distance of each of these points to $C^c$ is
$\geq n/2$ and yet, every $1$-chain joining them must pass through
$C_n\cap\{x_2=0\}$ whose points are at distance $\leq 2\log n$
from $C^c$. But this is impossible for a ball in a (M)-space.
Indeed, in a ball $B=B(o,R)$ with $R\geq N$, if a point $x$ is at
distance $cN$ from the boundary, then the points belonging to a
ball centered in $x$ and of radius $cN/2$ are at distance at least
$cN/2$ from the boundary of $B$. But this ball intersects the ball
centered in $o$ and of radius $R-cN/2$. Moreover, by property (M),
there exists a $1$-chain joining $x$ to $o$ and staying in
$B(o,R-cN/2)$, so at a distance of the order of $N$ from boundary
of $B$. \epr

\subsection{The instability of (IB) under quasi-isometry between graphs of polynomial growth}\label{S}

\begin{thm}\label{t1}
Let $d$ be an integer $\geq 2$. There exists two graphs $X$ and
$X'$ quasi-isometric, of polynomial growth of degree $d$ and with
bounded or unbounded profile, such that $X$ satisfies (IB) but not
$X'$.
\end{thm}

Like in the examples of the two previous sections, we will build a
graph $X$ removing some edges from $\Z^d$: for $n\in \N$, let
$A_n$ be the ball of radius $n$ whose center belongs to the axis
$\Z.e_1$ in such a chain that $A_{n+1}$ is at distance $2^n$ from
$A_n$. We then remove all the edges of the boundary of $A_n$ but
those belonging to the line $\Z.e_1$. We write $A$ for the union
of $A_n$. The graph $X'$ is obtained from $X$ by taking its image
by the linear map fixing the first coordinate and acting on the
orthogonal as an homothetic transformation of ratio $4$ (it is
clear that it is a quasi-isometry). More precisely, we replace
each edge of $X$ parallel to the first axis, by a chain of length
$2$ also parallel to the first axis. Write $A'$ for the image of
$A$.

\begin{rem}
In the previous example, the profile is bounded. Nevertheless, one
can slightly modify the construction in order to get an unbounded
profile: for instance, removing only edges of the boundary of
$A_n$ at distance $\geq \log n$ from the axis $\Z.e_1$ (instead of
those which are outside of this axis).
\end{rem}

\begin{clai}
The graphs $X$ and $X'$ have polynomial growth of degree $d$.
\end{clai}
As these graphs are subgraphs of $\Z^d$, their volume growths are
less than the one of $\Z^d$. The converse inequality will follow
from the fact that in $X'$, the profile restricted to balls is
strong and from the fact that $X$ and $X'$ are quasi-isometric.
\epr

\

\begin{clai}
In $X$, the balls are asymptotically isoperimetric.
\end{clai}
\bpr It is clear by construction that the $A_n$ are balls and that
their boundaries have bounded volume. \epr

\

\begin{clai}
In $X'$, the profile restricted to balls is strong:
$I^{\downarrow}_{\mathcal{B}}(t)\approx t^{\frac{d-1}{d}}$. In
particular, $X'$ is not (IB).
\end{clai}
\bpr Remark that Lemma \ref{l1} stays true in this context. Let
$B=B(x,r)$ be a ball of the graph $X'$. We have to show that there
exists a constant $c>0$ such that
$$|\partial B|\geq c|B|^{\frac{d-1}{d}}.$$

According to Lemma \ref{l1}, we can assume that $B\subset A$.
Thus, there exists $n_0$ such that $B\subset A_{n_0}$.

Let us embed $\Z^d$ into $\R^d$. Let us replace the discrete
polyhedron $A_n$ and $B$ by their convex hulls $\tilde{A_n}$ and
$\tilde{B}$ in $\R^d$. Let $\tilde{X}$ be the space obtained
removing from $\R^d$ (Euclidian) the points of the Euclidian
boundary of $\tilde{A_n}$ (for all $n$) but the two ones belonging
to the axis $\R.e_1$ (resp. those at distance $\leq \log n$ of the
axe) for the case of bounded profile (resp. for the case of
unbounded profile). Let us equip $\tilde{X}$ --seen as a subset of
$\R^d$-- with Lebesgue measure and with the geodesic metric
$d(x,y)=\inf_{\gamma} l(\gamma)$ with $\gamma$ taking values in
the set of arcs joining $x$ to $y$ in $\tilde{X}$, $l(\gamma)$
being the Euclidian length of $\gamma$.

The embedding $j$ of $X$ into $\tilde{X}$ we obtain like this is
clearly a large-scale equivalence.

For simplicity, we will write $|A|$ for the (Lebesgue) measure of
a subset $A$ of $\tilde{X}$. On the other hand, note that
$\partial_{10}\tilde{B}$ contains $[j(B(x,r))]_1\setminus
[j(B(x,r-2))]_1$, which by Proposition \ref{prop1} has same
measure (up to multiplicative constant) as $\partial B$. The same
holds for $\tilde{B}$ and $B$. Moreover, since $\tilde{B}$ and
$A_{n_0}$ are convex polyhedra, it is clear that the $10$-boundary
of $\tilde{B}$ has same measure (up to multiplicative constants)
as its Euclidian boundary (whose measure is the limit when
$h\rightarrow 0$ of $|\partial_h \tilde{B}|/h$). Write:
$$|\partial_{eucl}\tilde{B}|=\lim_{h\rightarrow 0}|\partial_h \tilde{B}|/h$$

Consequently, it is enough to show that there exists $c>0$ such
that
$$|\partial_{eucl}\tilde{B}|\geq c|\tilde{B}|^{\frac{d-1}{d}}$$
Note that by homogeneity, the quantity
$$Q=\frac{1}{r^{d-1}}|\partial_{eucl}\tilde{B}|$$ only depends on the
ratio $n/r$. Fix $n=n_0$. For $r$ small enough (let us say $\leq
r_c$ for some $r_c>0$), $\tilde{B}$ never meets two parallel
faces: $Q$ stays larger than a constant $>0$ (i.e. profile of a
$1/2^{d-1}$'th of space of $\R^d$). By compactness, it follows
that $Q$ reaches its minimum when $x$ and $r$ vary under the
conditions: $r_c\leq r\leq n_0/2$. On the other hand, as
$\tilde{B}$ is strictly included in $\tilde{A_{n_0}}$, this
minimum has to be $>0$. The ratio $Q$ is therefore larger than a
constant $c'>0$. finally, there is a constant $c>0$ such that:
$$|\partial_{eucl}\tilde{B}|\geq c' r^{d-1}\geq c|B|^{\frac{d-1}{d}}.$$
So we are done.\epr

\

\section{Asymptotic isoperimetry of connected
subsets.\\}\label{conn}

Recall that we say that a subset $A$ of a metric space is
connected if there does not exist a nontrivial partition
$A=A_1\sqcup A_2$ with $d(A_1,A_2)\geq 10$.

Let $(X,d,\mu)$ be a metric measure space. Write $\mathcal{C}$ for
the set of connected subsets of finite measure of $X$.

Set $\partial A=\partial_{1}A$ and assume that $X$ is uniformly
$1/2$-connected (see section \ref{loc}).

\begin{thm}

\

(i) Let $X$ be such that the measure of balls of radius $1/2$ are
minorated by $a>0$. Suppose that $I(t)=o(t)$. Then there exists a
positive and increasing sequence $(t_i)$ tending to infinity such
that $I^{\downarrow}_{\mathcal{C}}(t_i)= I(t_i)$.

(ii) Assume that $X$ is a doubling (M)-space and has a strong
profile. Then $I_c\approx I$.

(iii) Let $d$ be an integer $\geq 2$. there exists a graph $X$ of
polynomial growth of degree $d$ and a increasing sequence of
integers $(N_n)$ such that
$I^{\downarrow}_{\mathcal{C}}(N_n)=o(I(N_n))$.
\end{thm}
\bpr

Note that (ii) follows from Corollary \ref{p3} and from the fact
that property (M) implies that balls are connected.

Let us show the first assertion of the theorem. Suppose that there
exists $T>0$ such that $\forall t\geq T$,
$I(t)<I^{\downarrow}_{\mathcal{C}}(t)$. We will show that it
implies that:
\begin{equation}\label{equat1}
I(t)\geq a\frac{t}{T}.
\end{equation}

Write $t_m$ for the upper bound of the set of $t$ such that
$\forall s\leq t$, one has: $I(s)\geq a\frac{s}{T}$. Since $I$ is
nondecreasing, if $t_m$ is finite, then it is a maximum.

Remark that $t_m\geq T$ since the boundary of every nonempty
subset of $X$ contains a ball of radius $1/2$ (see Proposition
\ref{prop1}) and therefore has measure $\geq a$.

Suppose by absurd that $t_m$ is finite. By definition of $t_m$,
for all $s>t_m$ there exists a subset $A$ such that:
$$\mu(A)\geq s$$ and:
$$\mu(\partial A)<as/T.$$
Moreover, since $t_m\geq T$, we can suppose that:
$$\mu(\partial A)<I^{\downarrow}_{\mathcal{C}}(s)$$
(in particular, $A$ is not connected).

It follows that there exists a smallest positive integer $k$ such
that there exist $t_m\leq s\leq t_m+T/2$ and a subset $A$ of
measure $\geq s$, with $k$ connected components and whose boundary
has measure $<\min\{I^{\downarrow}_{\mathcal{C}}(s), sa/T\}$. Let
$A$ be such a subset. Note that $k\geq 2$. Thus, we have:

$$A=A_1\sqcup A_2$$
with $d(A_1,A_2)\geq 10$.

Since $k$ is minimal, one has, for $i=1,2$:
$$\mu(A_i)<t_m.$$
Indeed, if for instance, one had $\mu(A_1)\geq t_m$, then since
the boundary of $A_2$ has measure $\geq a$, one would have:
$$\mu(\partial A_1)\leq (t_m+T/2)\frac{a}{T}-a=\frac{t_ma}{T}-a/2<\frac{t_ma}{T}.$$
Therefore, as $I^{\downarrow}_{\mathcal{C}}(t_m)\geq I(t_m)\geq
t_ma/T$, one would also have:
$$\mu(\partial A_1)<I^{\downarrow}_{\mathcal{C}}(t_m).$$
But then, by minimality of $k$, $A_1$ should have at least $k$
connected components, which is absurd since it has strictly less
components than $A$.

But, by definition of $t_m$, this implies that:
\begin{eqnarray*}
\mu(\partial A) & = & \mu(\partial A_1)+\mu(\partial A_2)\\
                &\geq &\frac{\mu(A_1)a}{T}+\frac{\mu(A_2)a}{T}\\
                & = & \frac{\mu(A)a}{T}
\end{eqnarray*}
which is absurd. \epr

\

In order to show the second assertion of the theorem, we proceed
as in the previous sections: we start from the graph $\Z^d$, and
then we remove some edges. Let us consider the following family of
cubes $(C_n^m)_{0\leq m\leq n-1, n\in \N^*}$ of $\Z^d$: the
$C_n^m$ are Euclidian cubes of edges' length $2^{2^n}$ whose
centers are disposed along the axis $\Z.e_1$ as follows:
$C_n^{m+1}$ is the image of $C_n^m$ by the translation of vector
$n2^{2^n}.e_1$ and $C_n^{n-1}$ and $C_{n+1}^{1}$ are at distance
$(n+1)2^{2^{(n+1)}}$ to one another. To build the graph $X$, we
remove all the edges joining $C_n^m$ to the rest of the graph but
those which have a vertex belonging to the Euclidian cube $c_n^m$
of dimension $d-1$ of the boundary of $C_n^m$, of volume $2^{n^2}$
and centered in one of the two intersection points of $C_n^m$ with
the axis $\Z.e_1$. Write $C$ for the union of cubes $C_n^m$.

\begin{clai}
The growth in $X$ is polynomial of degree $d$.
\end{clai}

\bpr Let $B=B(x,r)$ be a ball. Let us prove that $|B|\geq
2^{-100d}r^d$. If the center of $B$ doesn't belong to any $C_n^m$,
it is clear. Suppose therefore that $x\in C_{n_0}^{m_0}$ for
integers $n_0$ and $m_0<n_0$. Write $D_{n_0}$ for the diameter of
$C_{n_0}^{m_0}$. If $r\geq 3D_{n_0}$, then  $B$ contains
$B(y,r/2)$ with $y$ belonging to no $C_n^m$. So we are brought
back to the previous case. In the other case, the conclusion
follows from the following trivial fact: in $\Z^d$, if $r\leq n$,
the volume of the intersection of a cube $C$ of edges' length
equal to $n$ with a ball of radius $r\leq n$ and of center $x\in
C$ is $\geq 2^{-d}|B(x,r)|\geq 2^{-10d}r^d$. Indeed, the worst
case is when $x$ is a corner of the cube. \epr

\

\begin{clai}
Take $N_n= n2^{2^n}$. Then
$I(N_n)=o(I^{\downarrow}_{\mathcal{C}}(N_n))$.
\end{clai}

\bpr Let us consider the set $C_n=\cup_m C_n^m$. Its volume is
equal to $N_n$ and its boundary has volume equal $n2^{n^2}$. On
the other hand, let $n_1$ be an integer and let $P$ be a connected
subset of volume $\geq N_{n_1}$. We want to show that $|\partial
P|\geq c2^{(n_1+1)^2}$, for a constant $c>0$, which is clearly
enough to conclude.

Thanks to the following lemma, the only remaining case to consider
is when $P$ meets a cube $C_n^m$. But, because of the large
distance between two such cubes, we can assume that $P$ meets only
one of these cubes, say $C_{n_0}^{m_0}$.
\begin{lem}
The profile of the graph $C^c$ is strong (i.e. $\approx
t^{\frac{d-1}{d}}$).
\end{lem}
(same demonstration as for Lemma \ref{l1})

\

If $|P\cap C^c|\geq|P|/2$, then the lemma applied to $P\cap C^c$
allows to conclude. Suppose therefore that $|P\cap C|\geq|P|/2$.
This implies in particular that $n_0\geq n_1+1$. We then remark
that $|\partial (P\cap C_{n_0}^{m_0})|\leq |\partial P|$. Indeed,
let $\pi$ be the orthogonal projection onto the hyperplane
containing $c_{n_0}^{m_0}$, then every point of $c_{n_0}^{m_0}\cap
P$ admits un antecedent by $\pi$ belonging to the boundary of $P$.
So we can assume that $P\subset C_{n_0}^{m_0}$. If $|P|\leq
3/2|C_{n_0}^{m_0}|$, then there exists $c>0$ such that:
\begin{equation}\label{lasteq}
|\partial P|\geq c|P|^{\frac{d-1}{d}}
\end{equation}
(isoperimetry in the full Euclidian cube: see \cite{Ros}).
Otherwise, assume that $|P|\geq 3/2|C_{n_0}^{m_0}|$ and write
$Q=C_{n_0}^{m_0}\setminus P$.

\begin{itemize}
\item If the volume of $Q$ is $\geq D_{n_0}/2$ where $D_{n_0}$ is
the diameter of $C_{n_0}^{m_0}$, then \ref{lasteq} applied to $Q$
implies that:
$$|\partial Q|\geq c2^{(d-1)2^{n_0}/d}\geq c2^{2^{n_0-1}}\geq
c2^{n_0^2}=c2^{(n_1+1)^2}.$$ But, the boundary of $Q$ is---up to
points belonging to $c_{n_0}^{m_0}$ (whose cardinal is negligible
compared to $c2^{2^{n_0-1}}$)---equal to the boundary volume of
$P$. So we are done.

\item If $|Q|\leq D_{n_0}/2$, then every point of $c_{n_0}^{m_0}$
has preimages in $\partial P$ by the projector $\pi$. But
$|c_{n_0}^{m_0}|= 2^{n_0^2}=2^{(n_1+1)^2}$, which ends the
demonstration.
\end{itemize}
\epr

\

\bibliographystyle{amsplain}

\bigskip
\footnotesize
\noindent Romain Tessera\\
Université de Cergy-Pontoise\\
E-mail: \url{tessera@clipper.ens.fr}\\

\end{document}